\documentclass[12pt,reqno]{article}
\usepackage{amsmath,amsfonts,amssymb,amsthm,indentfirst}

\catcode`\@=11

 \def\AMSTeXfeatures{\Plainheads 
   \let\current@vert=\AMS@vert}

 \def\Plainheads{\sh@ftdiam=0.05em
   \getlabeldims
   \let\vshaftfill=\plnvsolidfill
   \let\hshaftfill=\plnhsolidfill
   \let\th@rhead=\plnrhead
   \let\th@lhead=\plnlhead
   \let\th@dnhead=\plndnhead
   \let\th@uphead=\plnuphead}
 
 \def\glet{\global\let}

 \def\LaTeXfeatures{\catcode`\@=11
   \ifx\@clnwd\undefined \nol@g
      \input ltxcode.tex \dol@g \fi
   \ltxheads \let\current@vert=\new@vert
   \providelto \catcode`\@=\active}

 \def\nol@g{\def\wlog{\edef\garbage}}
 \def\dol@g{\let\wlog=\wl@g} \let\wl@g=\wlog
 \nol@g 

 \newbox\ltobox
 \def\providelto{{\setbox\z@=
   \hbox{$\to$}\minharrlen=\wd\z@
   \global\setbox\ltobox=\hbox{$\activeat>>>$}}
   \def\lto{\mathrel{\copy\ltobox}}}

 \def\ltxheads{\sh@ftdiam=\@wholewidth
   \getlabeldims
   \let\vshaftfill= \ltxvsolidfill
   \let\hshaftfill=\ltxhsolidfill
   \let\th@rhead=\ltxrhead
   \let\th@lhead=\ltxlhead
   \let\th@dnhead=\ltxdnhead
   \let\th@uphead=\ltxuphead}
 {\catcode`\@=\active
   \gdef@#1{\csname #1\string@at\endcsname}
   \glet\activeat=@}
 \def\def@#1{\expandafter\def\csname #1@at\endcsname}

 \def@>#1>#2>{\@rrow R{#1}{#2}}
 \def@<#1<#2<{\@rrow L{#1}{#2}}
 \def@ V#1V#2V{\@rrow V{#1}{#2}}
 \def@ A#1A#2A{\@rrow A{#1}{#2}}
 \def@/#1/#2/#3/{\@rrow{#1}{#2}{#3}}
 \def@.{\ifodd\row\ifmmode\noharrow
     \else\leavevmode.\spacefactor3000 \fi
   \else\novarrow\fi}
 \def@={\ifodd\row\harrow\hequalfill{}{}%
   \else\varrow\vequalfill{}{}\fi}
 \def@:#1{\ifx=#1\harrow\deffill{}{}%
   \else\leavevmode\null:#1\fi}
 \def@|{\current@vert}
  \def\AMS@vert{\varrow\vequalfill{}{}}
  \def\new@vert#1|#2|{\ifodd\row
   \let\nextarrow\vertexvarrow
   \else\let\nextarrow\varrow\fi
   \nextarrow\vshaftfill{#1}{#2}}
 \def@-{\ifmmode\let\next\hl@ne
   \else\let\next\AMSatdash \fi \next}
  \def\hl@ne#1-#2-{\harrow\hshaftfill{#1}{#2}}
  \def\AMSatdash{\let\next\relax\leavevmode
    \def\next@{\ifx\next-%
      \def\next-{\futurelet\next\nextii@}%
     \else\def\next{\hbox{-}}\fi\next}%
    \def\nextii@{\ifx\next-\def\next-{\hbox{---}}%
      \else\def\next{\hbox{--}}\fi\next}%
    \futurelet\next\next@}
 \def@(#1){\tweenarrows{#1}}
 \def@[#1]{\setsp@n#1\relax\activeat}
 \def\fiberbox{\hbox{$\vcenter{\hr@le\hbox{\vr@le
   \kern1ex\vbox{\kern1.2ex}\vr@le}\hr@le}$}}
  \def\hr@le{\hrule height \sh@ftdiam}
  \def\vr@le{\vrule width \sh@ftdiam}
 \def@+#1+#2+#3+{\ifodd\row \harrow{#1}{#2}{#3}%
   \else \varrow{#1}{#2}{#3}\fi}

 \def\Dnarrfill{\vequalfill\Dnhe@d}
 \def\Uparrfill{\Uphe@d\vequalfill}
 
 \def\ontofill{\rtarrfill\kern-0.3em 
   \th@rhead\kern 0.3em} 

 \def\rtarrfill{\hshaftfill\th@rhead}
 \def\ltarrfill{\th@lhead\hshaftfill}
 \def\dnarrfill{\vshaftfill\th@dnhead}
 \def\uparrfill{\th@uphead\vshaftfill}
 \def\hequalfill{\plnhfill=}
 \def\deffill{:\plnhfill=}
 \def\plnvextfill#1{\setbox\z@
   \hbox{\the\textfont3 #1}%
   \dimen@=\dp\z@\advance\dimen@\ht\z@
   \copy\z@ \kern-\dimen@ 
   \cleaders\copy\z@ \vfill
   \kern-\dimen@ 
   \box\z@}
 \def\plnhfill#1{$\m@th\mkern-1.5mu\mathord#1\mkern-6mu
    \cleaders\hbox{$\mkern-2mu\mathord#1\mkern-2mu$}\hfill
    \mkern-6mu\mathord#1\mkern-1.5mu$}
 \def\vequalfill{\plnvextfill{\char'167}}
 \def\plnvsolidfill{\plnvextfill{\char'077}}
 \def\plnhsolidfill{\plnhfill-}
 \def\ltxhsolidfill{\leaders\hrule height\topofshaft depth\botofshaft
   \hfill}
 \def\ltxvsolidfill{\leaders\vrule width\sh@ftdiam\vfill}
 \def\hdashfill{\hd@sh\wd@sh
   \xleaders \hbox{\wd@sh\hd@sh\wd@sh}\hfill
   \wd@sh\hd@sh}
 \def\vdashfill{\vd@sh\wd@sh
   \xleaders \vbox{\wd@sh\vd@sh\wd@sh}\vfill
   \wd@sh\vd@sh}
 \def\dashed{\ifinmeasureCD\else
    \ifodd\row\option{\let\hshaftfill=\hdashfill}%
   \else\option{\let\vshaftfill=\vdashfill}\fi\fi}

 \newdimen\CDstrutht  \newdimen\CDstrutdp
 \newdimen\CDstrutlen \CDstrutlen=\CDstrutht
   \advance\CDstrutlen by \CDstrutdp

 \def\CDstrut{\vrule
   height \ifnum\row=1 \z@\else\CDstrutht \fi
   depth \ifnum\row=\numrows \z@ \else\CDstrutdp \fi
   width\z@}

 \newdimen\CDarrsurr \CDarrsurr=0.375em
 \newdimen\CDdashlen
 \newdimen\CDvarrlen \CDvarrlen=1.5\baselineskip
 \newdimen\minharrlen 
  \setbox\z@\hbox{$\longrightarrow$} \minharrlen=\wd\z@
 \newdimen\minCDharrlen \minCDharrlen=2.5em 
\newdimen \minc@lwd
\def\findminc@lwd{\minc@lwd=2\CDarrsurr
  \advance\minc@lwd\minCDharrlen}

 \newdimen\sh@ftdiam

 \newdimen\labelsurr \labelsurr=1.25 em

\newcount\sp@ncnt \sp@ncnt=\@ne
\newcount\sp@ncnt@ \sp@ncnt@=\@ne
\newdimen\@rrwd \newdimen\@rrdp

 \def\adjustbot#1{\option{\advance\@rrdp#1\relax}}
 
\def\pushvertex#1{\global\p@shlen#1\relax
   \global\let\maybepush=\dopush}

 \newdimen\p@shlen \p@shlen=\z@

 \let\maybepush=\relax
 \def\dopush{\ifinmeasureCD 
   \advance\locdimen by -\p@shlen 
   \else\advance \@rrwd by -\p@shlen \fi 
   \global\let\maybepush=\relax \global\p@shlen=\z@\relax}

 \def\span@ne{\global\sp@ncnt=\@ne\relax}
 \def\setsp@n#1#2{\global\sp@ncnt=#1\relax
   \ifx\relax#2\relax\else\global\sp@ncnt@=#2\relax\fi}

 \def\plnrhead{\llap{$\rightarrow\mkern-1.5mu$}}
 \def\plnlhead{\rlap{$\mkern-1.5mu\leftarrow$}}

 \def\clap#1{\hbox to \z@{\hss #1\hss}}

 \def\plndnhead{\hbox{\the\textfont3 \char'171}}
 \def\plnuphead{\hbox{\the\textfont3 \char'170}}
 \def\Dnhe@d{\hbox{\the\textfont3 \char'177}}
 \def\Uphe@d{\hbox{\the\textfont3 \char'176}}

 \def\ltxrhead{\raise\@xisheight
   \llap{\smash{\@linefnt\@getrarrow(1,0)}}}
 \def\ltxlhead{\raise\@xisheight
   \rlap{\@linefnt\@getlarrow(-1,0)}}
 \def\ltxuphead{\setbox\z@=\rlap{%
   \kern\@halfwidth\@linefnt\char'66}%
   \copy\z@\kern-\ht\z@}
 \def\ltxdnhead{\setbox\z@=\rlap{%
   \kern\@halfwidth\@linefnt\char'77}%
   \ht\z@=\z@\box\z@}

 \def\wd@sh{\kern0.5\CDdashlen}
 \def\hd@sh{\vrule height\topofshaft depth\botofshaft
    width\CDdashlen}
 \def\vd@sh{\hrule height\CDdashlen
   depth\z@ width\sh@ftdiam}

\def\xylist{14{3434}13{2414}12{1723}%
  23{1413}34{1153}11{0867}43{0707}%
  32{0580}21{0414}31{0291}41{0}}
\newcount\tgtcnt@
\def\find@xyargs{\dimen@=\@rrdp
  \advance\dimen@ by \CDstrutlen
  \tgtcnt@=\dimen@ \dimen@=\@rrwd 
  \divide\dimen@ by \@m 
  \divide \tgtcnt@ by \dimen@ 
  \expandafter\testxy\xylist\relax
  \unitlength=\@xarg\@rrdp
  \divide\unitlength by\@yarg\relax}
\def\testxy#1#2#3{\ifnum\tgtcnt@>#3
    \@xarg=#1\relax \@yarg=#2\relax
    \let\next=\ignorerest
  \else\let\next\testxy\fi\next}
\def\ignorerest#1\relax{\relax}

\let\scalefactor=\@ne
\def\SWarrow{\find@xyargs\vector
  (-\@xarg,-\@yarg)\scalefactor\hskip-\wd\@linechar}
\def\NWarrow{\find@xyargs\vector
  (-\@xarg,\@yarg)\scalefactor\hskip-\wd\@linechar}
\def\NEarrow{\find@xyargs\vector
  (\@xarg,\@yarg)\scalefactor}
\def\SEarrow{\find@xyargs\vector
  (\@xarg,-\@yarg)\scalefactor}
\def\rightupline{\find@xyargs\@linelen=\scalefactor
     \unitlength\@sline}
\def\rightdownline{\find@xyargs\@yarg=-\@yarg\relax
     \@linelen=\scalefactor\unitlength\@sline}

\def\Sim{\ifodd\row\setbox\z@=\hbox{$\sim$}\dimen@=\ht\z@
 \advance\dimen@ by -\@xisheight
  \vbox{\box\z@\kern-\@xisheight\kern\dimen@}%
  \else\hbox{$\wr$}\fi}

\def\harrow#1#2#3{\inmeasureCDtrue\findminarrwd
  {#2}{#3}{\sp@ncnt\minharrlen}\inmeasureCDfalse\span@ne
  \mathrel{\hbox{\options\hplace{#1}\ulabel{#2}\dlabel{#3}}}}

\def\noharrow{\harrow\hfill{}{}}
\def\vertexvarrow#1#2#3{\findarrdp \@rrwd=\z@ \setsp@n\@ne\@ne
  \vbox to \z@{\kern-1.2\CDstrutht
  \rlap{\options\vplace{#1}\llabel{#2}\rlabel{#3}}\vss}}

\newif\ifinmeasureCD
\def\measurelabel#1{\setbox\z@
  \hbox{$\scriptstyle#1\kern\labelsurr$}%
  \ifdim\wd\z@>\@rrwd \@rrwd=\wd\z@\fi}
\def\findminarrwd#1#2#3{\@rrwd=#3\relax
   \measurelabel{#1}\measurelabel{#2}}
\def\findCDarrwd#1#2{\@rrwd=\minCDharrlen
   \measurelabel{#1}\measurelabel{#2}%
  }

\newcount\row \row=\@ne \newcount\col \col=\@ne 
 \newcount\numrows 
\numrows=\@ne
 \newcount\numcols
\newcount\arrspan \newdimen\vrtxhalfwd  \newbox\tempbox

\def\DANABUG{\advance\col by \@ne
 \@rrwd=\minCDharrlen
  \advance\@rrwd by \vrtxhalfwd
  \advance\@rrwd by \CDarrsurr
  \ifnum\col>\numcols \numcols=\col
     \newlocdimen{col\the\col}\locdimen=\@rrwd 
  \else \ifdim\@rrwd>\c@l \c@l=\@rrwd\fi\fi}

\def\drop#1\\{
  \findvrtxhalfsum\DANABUG\advance\row by 2 \measureinit}

\def\measureinit{\col=\@ne \vrtxhalfwd=-\CDarrsurr\arrspan=\@ne\@rrwd=\z@
   \setbox\tempbox=\hbox\bgroup$}
\def\measure{
  \let\harrow\measureCDarrow
  \let\CDCR=\measureCR 
   \findminc@lwd 
  \inmeasureCDtrue
  \row=\@ne \numcols=\z@ \measureinit}

\def\endmeasure{\findvrtxhalfsum\DANABUG
  \numrows=\row 
  \inmeasureCDfalse}

\def\newlocdimen#1{\advance\dimenc@unt by \@ne
  \ifnum\dimenc@unt<\insc@unt
     \else\errmessage{No room for the CD}\fi
  \dimendef\locdimen=\dimenc@unt
  \expandafter\dimendef\csname#1\endcsname=\dimenc@unt}

 \def\r@wc@l{\csname row\the\row col\the\col\endcsname}
 \def\c@l{\csname col\the\col\endcsname}

 \def\findvrtxhalfsum{$\egroup
  \newlocdimen{row\the\row col\the\col}
  \locdimen=\vrtxhalfwd 
  \vrtxhalfwd=0.5\wd\tempbox 
  \advance\vrtxhalfwd by \CDarrsurr
  \advance\locdimen by \vrtxhalfwd 
  \advance\@rrwd by \locdimen 
  \maybepush
  \divide\@rrwd by \arrspan\relax
  \ifdim\@rrwd<\minc@lwd
    \ifnum\col>\@ne \@rrwd=\minc@lwd\fi \fi
  \loop 
    \ifnum\col>\numcols \numcols=\col
       \newlocdimen{col\the\col}
       \locdimen=\@rrwd 
    \else \ifdim\@rrwd>\c@l \c@l=\@rrwd\fi \fi
   \ifnum\arrspan>\@ne
      \advance\arrspan by -1 \advance\col by \@ne
  \repeat }

 \def\measureCDarrow#1#2#3{\findvrtxhalfsum
   \arrspan=\sp@ncnt\relax\global\sp@ncnt=1\relax
   \advance\col by \@ne
   \findCDarrwd{#2}{#3}%
   \setbox\tempbox=\hbox\bgroup$}

 \newcount\dr@tn \dr@tn=\z@
 \def\locate#1:#2{\ifinmeasureCD\else
   \count@=-#1
   \multiply\count@ by 2
   \advance\count@ by #2
   \dimen@=\count@\@rrwd
   \ifnum\dr@tn=\@ne\relax \else\dimen@=-\dimen@ \fi
   \dimen@i=\@rrdp
   \ifnum\dr@tn>\z@\advance\dimen@i by \CDstrutlen \fi
   \dimen@i=\count@\dimen@i
   \count@=#2 \multiply\count@ by 2
   \divide\dimen@ by \count@
   \divide\dimen@i by \count@
   \lift\dimen@i\nudge\dimen@\fi}

\def\betweenCDrows{\advance\row by \@ne \col=\@ne
\options}

\def\hbegin{\hbox\bgroup\kern\c@l \kern-\r@wc@l$}
\def\hend{$\glet\maybepush\relax \CDstrut\egroup}
\def\vbegin{\setbox\tempbox=\hbox\bgroup$}
\def\vend{$\egroup\ht\tempbox=\z@\dp\tempbox\CDvarrlen
  \box\tempbox}
\def\setCD{\let\harrow=\setCDarrow
  \let\CDCR=\setCR 
  \row=\@ne \col=\@ne \hbegin}
\let\endsetCD=\hend 

\def\findarrwd{\@rrwd=\z@ \count@=\col \advance\count@ by\sp@ncnt
  \loop\ifnum\count@>\col \advance\count@ by -1
      \advance\@rrwd by\csname col\the\count@\endcsname\repeat}
\def\setCDarrow#1#2#3{\kern\CDarrsurr\advance\col by \@ne
  \findarrwd \advance\@rrwd by -\r@wc@l  
  \@rrdp=\z@ 
  \maybepush
  \advance\col by -\@ne \advance\col by \sp@ncnt \span@ne
  \hbox to \@rrwd{\options
   \@rrwd=\scalefactor\@rrwd\hss
   \hplace{#1}\ulabel{#2}\dlabel{#3}\hss}%
   \kern\CDarrsurr}

\newdimen\labspacei 
\newdimen\labspaceii 

\newdimen\@xisheight
  \@xisheight=\the\fontdimen22\textfont2
\newdimen\labelskip
  \labelskip=\the\fontdimen10\textfont3 
\newdimen\topofshaft
\newdimen\botofshaft
\newdimen\botofulabel
\newdimen\topofdlabel
\def\getlabeldims{
  \topofshaft=0.5\sh@ftdiam
  \botofshaft=\topofshaft
  \advance\topofshaft by \@xisheight  
  \advance\botofshaft by -\@xisheight  
  \botofulabel=\topofshaft
  \advance\botofulabel by \labelskip
  \topofdlabel=\botofshaft
  \advance\topofdlabel by \labelskip}

\def\ulabel{\ifnum\row=\@ne\let\next\ulabeli
   \else\let\next\ulabellap\fi\next}
\def\ulabeli#1{\vbox{
  \clap{\kern-\@rrwd$\scriptstyle#1$}%
  \kern\botofulabel}\maybeoffset}
\def\ulabellap#1{\vbox to \z@{\vss
  \clap{\kern-\@rrwd$\scriptstyle#1$}%
  \kern\botofulabel}\maybeoffset}
\def\dlabel{\ifnum\row=\numrows\let\next\dlabeli
   \else\let\next\dlabellap\fi\next}
\def\dlabeli#1{\vtop{\kern\topofdlabel
  \clap{\kern-\@rrwd$\scriptstyle#1$}%
  }\maybeoffset}
\def\dlabellap#1{\vbox to \z@{\kern\topofdlabel
  \clap{\kern-\@rrwd$\scriptstyle#1$}%
  \vss}\maybeoffset}
\def\rlabel#1{\vbox to \z@{\vss
  \rlap{\kern\labelskip$\scriptstyle#1$}%
  \vss\kern-\@rrdp}\maybeoffset}
\def\llabel#1{\vbox to \z@{\vss
  \llap{$\scriptstyle#1$\kern\labelskip}%
  \vss\kern-\@rrdp}\maybeoffset}
\def\swlabel#1{\vtop{\kern0.5\@rrdp
  \llap{$\scriptstyle#1$\kern\labelskip\kern-0.5\@rrwd}
  }\maybeoffset}
\def\nwlabel#1{\vbox{
  \llap{$\scriptstyle#1$\kern\labelskip\kern-0.5\@rrwd}%
  \kern-0.5\@rrdp}\maybeoffset}
\def\selabel#1{\vtop{\kern0.5\@rrdp
  \rlap{\kern0.5\@rrwd\kern\labelskip$\scriptstyle#1$}%
  }\maybeoffset}
\def\nelabel#1{\vbox{
  \rlap{\kern0.5\@rrwd\kern\labelskip$\scriptstyle#1$}%
  \kern-0.5\@rrdp}\maybeoffset}
\def\cplace#1{\vbox to \z@{\vss
  \clap{$#1$\kern-\@rrwd}%
  \kern-\@rrdp\vss}\maybeoffset}
\def\hplace#1{\hbox to \@rrwd{#1}\maybeoffset}
\def\vplace#1{\clap{\vbox to \z@{#1\kern-\@rrdp}}\maybeoffset}

\newdimen\nudgeamount \nudgeamount=\z@
\newdimen\liftamount \liftamount=\z@
\let\maybeoffset\relax
\newbox\offsetbox \newdimen\lastheight
\def\dooffset{
  \setbox\offsetbox=\lastbox \lastheight=\ht\offsetbox 
  \setbox\offsetbox=\vbox{\kern-\liftamount\box\offsetbox}%
  \ht\offsetbox=\lastheight
  \kern\nudgeamount\box\offsetbox\kern-\nudgeamount
  \global\nudgeamount=\z@ \global\liftamount=\z@
  \glet\maybeoffset=\relax}
\def\nudge#1{\ifinmeasureCD\else
  \global\advance\nudgeamount#1\relax
  \global\let\maybeoffset\dooffset\fi}
\def\lift#1{\ifinmeasureCD\else
  \global\advance\liftamount#1\relax
  \global\let\maybeoffset\dooffset\fi}

\def\findarrdp{\@rrdp=\CDvarrlen
  \ifnum\sp@ncnt@>1
    \advance\@rrdp by \CDstrutlen
    \multiply\@rrdp by \sp@ncnt@
    \advance\@rrdp by -\CDstrutlen \fi
 }

\def\varrow#1#2#3{\ifnum\sp@ncnt>\@ne 
     \sp@ncnt@=\sp@ncnt\relax\fi
  \findarrdp \@rrwd=\z@ 
  \kern\c@l
   \hbox to \z@{\options
   \@rrdp=\scalefactor\@rrdp
    \hss\vplace{#1}\llabel{#2}\rlabel{#3}\hss}%
  \global\advance\col by \@ne \setsp@n\@ne\@ne
  }

\def\novarrow{\varrow\vfill{}{}}

\def\tweenarrows#1{\findarrwd \findarrdp \setsp@n\@ne\@ne
  \rlap{\options\cplace{#1}}}

\def\usarrow #1#2#3{\dr@tn=\@ne
  \findarrwd \findarrdp \setsp@n\@ne\@ne 
  \rlap{\options\cplace{#1}\nwlabel{#2}\selabel{#3}}%
  \dr@tn=\z@}
\def\dsarrow #1#2#3{\dr@tn=\tw@
  \findarrwd \findarrdp \setsp@n\@ne\@ne 
  \rlap{\options\cplace{#1}\swlabel{#2}\nelabel{#3}}%
  \dr@tn=\z@}
 \def\@rrow#1{\csname #1@rrow\endcsname}
 \def\R@rrow{\harrow \rtarrfill}
 \def\L@rrow{\harrow \ltarrfill}
 \def\V@rrow{\varrow \dnarrfill}
 \def\A@rrow{\varrow \uparrfill}
 \def\SE@rrow{\dsarrow \SEarrow}
 \def\NW@rrow{\dsarrow \NWarrow}
 \def\SW@rrow{\usarrow \SWarrow}
 \def\NE@rrow{\usarrow \NEarrow}
 \def\DS@rrow{\dsarrow \dnslope}
 \def\US@rrow{\usarrow \upslope}
 \def\upslope{\find@xyargs
       \@linelen=\unitlength\@sline}
 \def\dnslope{\find@xyargs\@yarg=-\@yarg\relax
       \@linelen=\unitlength\@sline}

\newtoks\optionlist 
\optionlist={}
\let\options\relax
\def\dooptions{\the\optionlist\global\optionlist={}%
  \glet\options=\relax}
\def\option#1{\ifinmeasureCD\else
  \glet\options=\dooptions
  \global\optionlist=\expandafter{\the\optionlist\relax#1}\fi}
\def\wider#1{\ifinmeasureCD\else
   \option{\advance\@rrwd by #1}\fi}
\def\deeper#1{\ifinmeasureCD\else
   \option{\advance\@rrdp by #1}\fi}

{\def\\{\global\let\sptoken= }\\ }

\def\CR{\futurelet\nexttok\testCR}
\def\testCR{\ifx\nexttok\sptoken
   \let\next\eatspaceCR\else\let\next\CDCR\fi\next}
\def\eatspaceCR#1 {\CR}
\def\measureCR{\ifx\nexttok\endmeasure\let\nextCR\relax
    \else\let\nextCR\drop\fi\nextCR}
\def\setCR{\ifodd\row
  \ifx\nexttok\endsetCD\else\hend\betweenCDrows\vbegin\fi
  \else\vend\betweenCDrows\hbegin\fi}

\countdef\dimenc@unt=11
\def\CD#1\endCD{
   \begingroup\let\\=\CR
  \m@th\offinterlineskip
   \measure#1\endmeasure\null\,\vcenter{\setCD#1\endsetCD}\,
   \endgroup
    }

\ifx\@clnwd\undefined \nol@g\else\catcode`\ =14\relax\fi
 \font\@linefnt=line10 
 \newcount\@tempcnta
 \newcount\@tempcntb
 \newdimen\@tempdima
 \newdimen\@tempdimb
 \newdimen\@wholewidth
 \newdimen\@halfwidth
   \@wholewidth\fontdimen8\@linefnt \@halfwidth .5\@wholewidth
 \newdimen\unitlength
 \newcount\@xarg
 \newcount\@yarg
 \newcount\@yyarg
 \newbox\@linechar
 \newdimen\@linelen
 \newdimen\@clnwd
 \newdimen\@clnht
 \newif\if@negarg
 
 \def\@whilenoop#1{}

 \def\@whiledim#1\do #2{\ifdim #1\relax#2\@iwhiledim{#1\relax#2}\fi}

 \def\@iwhiledim#1{\ifdim #1\let\@nextwhile=\@iwhiledim 
         \else\let\@nextwhile=\@whilenoop\fi\@nextwhile{#1}}

 \def\@sline{\ifnum\@xarg< 0 \@negargtrue \@xarg -\@xarg \@yyarg -\@yarg
   \else \@negargfalse \@yyarg \@yarg \fi
 \ifnum \@yyarg >0 \@tempcnta\@yyarg \else \@tempcnta -\@yyarg \fi
 \ifnum\@tempcnta>6 \@badlinearg\@tempcnta0 \fi
 \ifnum\@xarg>6 \@badlinearg\@xarg 1 \fi
 \setbox\@linechar\hbox{\@linefnt\@getlinechar(\@xarg,\@yyarg)}%
 \ifnum \@yarg >0 \let\@upordown\raise \@clnht\z@
    \else\let\@upordown\lower \@clnht \ht\@linechar\fi
 \@clnwd=\wd\@linechar
 \if@negarg \hskip -\wd\@linechar \def\@tempa{\hskip -2\wd\@linechar}\else
      \let\@tempa\relax \fi
 \@whiledim \@clnwd <\@linelen \do
   {\@upordown\@clnht\copy\@linechar
    \@tempa
    \advance\@clnht \ht\@linechar
    \advance\@clnwd \wd\@linechar}%
 \advance\@clnht -\ht\@linechar
 \advance\@clnwd -\wd\@linechar
 \@tempdima\@linelen\advance\@tempdima -\@clnwd
 \@tempdimb\@tempdima\advance\@tempdimb -\wd\@linechar
 \if@negarg \hskip -\@tempdimb \else \hskip \@tempdimb \fi
 \multiply\@tempdima \@m
 \@tempcnta \@tempdima \@tempdima \wd\@linechar \divide\@tempcnta \@tempdima
 \@tempdima \ht\@linechar \multiply\@tempdima \@tempcnta
 \divide\@tempdima \@m
 \advance\@clnht \@tempdima
 \ifdim \@linelen <\wd\@linechar
    \hskip \wd\@linechar
   \else\@upordown\@clnht\copy\@linechar\fi}
 
 \def\@getlinechar(#1,#2){\@tempcnta#1\relax\multiply\@tempcnta 8
 \advance\@tempcnta -9 \ifnum #2>0 \advance\@tempcnta #2\relax\else
 \advance\@tempcnta -#2\relax\advance\@tempcnta 64 \fi
 \char\@tempcnta}
 
 \def\vector(#1,#2)#3{\@xarg #1\relax \@yarg #2\relax
 \@tempcnta \ifnum\@xarg<0 -\@xarg\else\@xarg\fi
 \ifnum\@tempcnta<5\relax
 \@linelen=#3\unitlength
 \ifnum\@xarg =0 \@vvector 
   \else \ifnum\@yarg =0 \@hvector \else \@svector\fi
 \fi
 \else\@badlinearg\fi}
 
 \def\@svector{\@sline
 \@tempcnta\@yarg \ifnum\@tempcnta <0 \@tempcnta=-\@tempcnta\fi
 \ifnum\@tempcnta <5
   \hskip -\wd\@linechar
   \@upordown\@clnht \hbox{\@linefnt  \if@negarg 
   \@getlarrow(\@xarg,\@yyarg) \else \@getrarrow(\@xarg,\@yyarg) \fi}%
 \else\@badlinearg\fi}
 
 \def\@getlarrow(#1,#2){\ifnum #2 =\z@ \@tempcnta='33\else
 \@tempcnta=#1\relax\multiply\@tempcnta \sixt@@n \advance\@tempcnta
 -9 \@tempcntb=#2\relax\multiply\@tempcntb \tw@
 \ifnum \@tempcntb >0 \advance\@tempcnta \@tempcntb\relax
 \else\advance\@tempcnta -\@tempcntb\advance\@tempcnta 64
 \fi\fi\char\@tempcnta}
 
 \def\@getrarrow(#1,#2){\@tempcntb=#2\relax
 \ifnum\@tempcntb < 0 \@tempcntb=-\@tempcntb\relax\fi
 \ifcase \@tempcntb\relax \@tempcnta='55 \or 
 \ifnum #1<3 \@tempcnta=#1\relax\multiply\@tempcnta
 24 \advance\@tempcnta -6 \else \ifnum #1=3 \@tempcnta=49
 \else\@tempcnta=58 \fi\fi\or 
 \ifnum #1<3 \@tempcnta=#1\relax\multiply\@tempcnta
 24 \advance\@tempcnta -3 \else \@tempcnta=51\fi\or 
 \@tempcnta=#1\relax\multiply\@tempcnta
 \sixt@@n \advance\@tempcnta -\tw@ \else
 \@tempcnta=#1\relax\multiply\@tempcnta
 \sixt@@n \advance\@tempcnta 7 \fi\ifnum #2<0 \advance\@tempcnta 64 \fi
 \char\@tempcnta}
\catcode`\ =10

\dol@g 
\catcode`\@=\active
\LaTeXfeatures

\newtheorem{thm}{Theorem}
\newtheorem*{prob*}{Problem}

\newtheorem{prop}{Proposition}[section]

\theoremstyle{definition}

\newtheorem*{rem*}{Remark}

\def\vp{\varphi}
\def\a{\alpha}

\def\cx{\mathcal{X}}

\numberwithin{equation}{section}
\date{}
\title{Geometrical equivalence, geometrical similarity\\ and
geometrical compatibility of algebras}
\author{Boris Plotkin}

\begin{document}

\maketitle
\begin{abstract}
In the paper the main attention is paid to conditions on algebras from
a given variety $\Theta$ which provide coincidence of their
algebraic geometries.
  The main part here
play the notions mentioned in the title of the paper.
\end{abstract}

\section{Introduction}  Let $\Theta$ be an arbitrary variety of algebras
 and $H$ an algebra from  $\Theta$. Then one can speak of
geometry in $\Theta$ over $H$.  We consider the classical
algebraic geometry with the  given field of coefficients $P$ as a
geometry, associated with the variety of all commutative and
associative algebras with a unit over $P$.  Denote this variety by
$Com-P$. The corresponding algebras $H$  are various extensions
$L$ of the field $P$.

For an arbitrary variety $\Theta$ we are looking for conditions
that provide the coincidence of geometries over two different
algebras $H_1$ and $H_2$ from $\Theta$.  We study this problem
from general positions, and consider also the cases of particular
varieties $\Theta$. In this paper we treat the following cases:
\begin{enumerate}
\item[1.] Classical case $Com-P$.  The problem is solved.

\item[2.] The case $\Theta = Ass-P$, all associative, not
necessarily commutative algebras over $P$.   The problem is open,
however, there exists a very clear conjecture.

\item[3.] Lie-$P$, all Lie algebras over $P$.   The problem is
solved.
\end{enumerate}

In all these cases the notions from the title play the crucial
role.

We will give all the necessary definitions and explain how to
understand that two geometries are the same.  It could be done in
different ways.

This paper is the full text of the talk, given at the 21st of
September, 2002, in St. Petersburg. Most of the proofs are
contained in the forthcoming paper [23].
 The presented paper is shorter and reflects the general picture
 in more clear and compact way. It also
contains some proofs which are absent in [23].
 Besides, we formulate a list of
problems, stimulated by the idea of coincidence of geometries.

In the theory we are working with the emphasis is made on
equations and their solutions in algebras from the given $\Theta$.
We study geometrical properties of algebras from $\Theta$ and
geometrical relations between algebras.

For special cases when $\Theta = Grp$ is a variety of all groups
and $\Theta=Grp-F$ is a variety of groups with the fixed free
group of constants, the corresponding geometry is related to
investigations of the Tarski's problem on elementary theory of a
free group.  There is a huge bibliography on this topic [9,10, 25
], etc..

For every variety $\Theta$ and algebra $H \in \Theta$ consider a
category $K_\Theta(H)$ of all algebraic sets over $H$.  Denote by
$\tilde K_\Theta(H)$ the category of algebraic varieties over $H$.
Algebraic variety is an algebraic set, considered up to
isomorphisms of the category $K_\Theta(H)$.  Hence, the category
$\tilde K_\Theta(H)$ is the skeleton of the category
$K_\Theta(H)$.

Both these categories are geometrical invariants of the algebra
$H$ and, to some extent, are responsible for the geometry in $H$.
The problem when the geometries over $H_1$ and $H_2$ are the same
is specified in the following two problems:

{\bf Problem 1}.  When are the categories $K_\Theta(H_1)$ and
$K_\Theta(H_2)$  isomorphic?

The second problem concerns isomorphism of categories of algebraic
varieties.  Recall  that in category theory it is proved that the
skeletons of two categories are isomorphic if and only if these
categories are equivalent.  Thus, we come to:

{\bf Problem 2}. When are the categories $K_\Theta (H_1)$ and
$K_\Theta(H_2)$  equivalent?

For every variety $\Theta$ consider also categories $K_\Theta$ and
$\tilde K_\Theta$.  Here the algebra $H$ is not fixed.  Both these
categories are geometrical invariants of the whole variety
$\Theta$.

{\bf Problem 3}.  When are $K_{\Theta_1}$ and $K_{\Theta_2}$
isomorphic or equivalent?

Here, the varieties $\Theta_1$ and $\Theta_2$ may be subvarieties
of a larger $\Theta$.

\section{Definitions}
\subsection{Varieties, prevarieties and quasivarieties of
algebras}

Recall that a variety of algebras is a class of algebras,
determined by a set of identities in some signature.  If $\cx$ is
an arbitrary class of algebras in some signature, then $Var (\cx)$
is a variety of algebras, determined by identities of the class
$\cx$. The variety $Var(\cx)$ is said to be generated by the class
$\cx$. The theorem
\[\tag{G. Birkhoff}
Var(\cx) = QSC(\cx)
\]
holds.  Here, $Q, S$ and $C$ are operators on classes of algebras,
where $C$ takes Cartesian products of algebras, $S$ takes
subalgebras and $Q$ takes homomorphic images.

In every variety of algebras $\Theta$ for every set of variables
$X$ there is the free algebra $W=W(X)$.  This is important for
logic and  geometry in $\Theta$.  Given $\Theta$, denote by
$\Theta^0$ the category, whose objects are free in $\Theta$
algebras $W=W(X)$ with finite $X$, and morphisms are homomorphisms
of algebras.

Along with varieties we consider prevarieties.  These are classes,
closed under operators $S$ and $C$.  If $\cx$ is an arbitrary
class of algebras, then the prevariety, generated by this class is
denoted by $pVar(\cx)$.  We have: $pVar(\cx) = SC(\cx)$. For every
class $\cx$ we consider also a local operator $L$ defined by the
rule: $H\in L\cx$ if every finitely-generated subalgebra in $H$
belongs to $\cx$.  We call a prevariety $\cx$ locally closed, if
it is closed in respect to the operator $L$. The locally closed
prevariety generated by a $\cx$ is equal to $LSC(\cx) $.

Consider further  quasiidentities and quasivarieties.
Quasiidentities are the formulas of the form
\[\tag{$\ast$} (\bigwedge_T (w\equiv w')) \Rightarrow w_0 \equiv w'_0 \; \;\; (w, w')
\in T
\]

Here, $T$ is a binary relation in $W=W(X)$, $X$ is finite  and all
$w, w',w_0,w_0'$ are elements of $W$. If the set $T$ is finite,
then $(\ast)$ defines  a usual quasiidentity.  In the general
case, we call the formulas  of the form $(\ast)$ infinitary
quasiidentities.

A quasivariety is the class of algebras defined  by a set of
quasiidentities. Let now $\cx$ be a subclass in $ \Theta$. Denote
by  $qVar(\cx)$ the quasivariety, determined by quasiidentities of
the class $\cx$. We have (see {A. Maltsev [13], Gretzer-Lakser
[7])
$$
 qVar (\cx) = SCC_{up} (\cx).
$$
Here, $C_{up}$ is the operator of taking of ultraproducts.  We
will use this formula in the sequel.



Denote by $\tilde qVar (\cx)$ the class defined by infinitary
quasiidentities of the class $\cx$.  The following theorem holds:
\begin{thm} (see [22])
\[
\tilde q Var (\cx) = LSC(\cx)
\]
for every $\cx \subset \Theta$.
\end{thm}

 Varieties and quasivarieties are
axiomatizable classes.  Classes $LSC(\cx)$ are not axiomatizable
in general.  They turn to be  axiomatizable in infinitary logic.
The following inclusions take place
\[
\tilde q Var (\cx) \subset q Var (\cx) \subset Var (\cx).
\]
  One of the central problems here is to find conditions providing
\[
\tilde q Var (\cx) = qVar (\cx).
\]

This problem was inspired by A. Maltsev and investigated by V.
Gorbunov [6]

\subsection{Affine spaces}

Fix a variety $\Theta$.  Take an algebra $H\in \Theta$ and a free
in $\Theta$ algebra $W=W(X)$ with finite $X$.  The set of
homomorphisms $Hom (W, H)$ we consider as an affine space of
points over $H$.  Points here are homomorphisms $\mu:W\to H$.  If
$X = \{ x_1, \cdots, x_n\}$, then we have a bijection
\[
\alpha_X :Hom(W, H) \to H^{(n)},
\]
and $\a_X(\mu) = (\mu(x_1),\cdots,\mu(x_n))$.  The point $\mu $ is
a root of the pair $(w, w')$, $w, w'\in W$, if $w^\mu=w'^\mu$,
which means also that $(w, w') \in Ker \mu$.  Here $Ker \mu$ is,
in general, a congruence of the algebra $W$.  Simultaneously,
$\mu$ is a solution of the equation $w = w'$.  We will identify
the pair $w, w'$ and the equation $w = w'$.

\subsection{Galois correspondence}

Let $T$ be a system of equations in $W$ and $A$ a set of points in
$Hom(W, H)$..  We have the following Galois correspondence
\[
\left\{
\begin{aligned}&T'_H =\{ \mu:W\to H\ | T\subset Ker\mu\}\\
&A'_W =\bigcap\limits_{\mu\in A} Ker \mu
\end{aligned} \right.
\]
The set $A$ of the form $A=T'$ for some $T$ we call a (closed)
algebraic set.  The congruence $T$ of the form $T=A'$ for some $A$
is an $H$-closed congruence.

It is easy to see that the congruence $T$ is $H$-closed if and
only if $W/T \in SC(H)$.

One can consider the closures $A^{''} = (A')'$ and $T_H^{''} =
(T'_H)'$.
\begin{prop}
The pair $(w_0, w'_0)$ belongs to  $T^{''}_H$ if and only if the
formula
\[
\Big(\bigwedge\limits_{((w, w')\in T} (w\equiv w')\Big)
\Rightarrow w_0 \equiv w'_0
\]
holds in $H$.
\end{prop}
\subsection{Categories}

We have defined the category $\Theta^0$.  Let us add to the
definition that all objects of $\Theta^0$ i.e., all finite $X$ are
subsets of an infinite universum $X^0$. Then $\Theta^0$ is a small
category.

Define further the category of affine spaces $K^0_\Theta(H)$.
Objects of this category are affine spaces
\[
Hom(W, H),\; \;  W \in 0b\; \Theta^0.
\]
The morphisms
\[
\tilde s: Hom(W(X), H) \to Hom (W(Y), H)
\]
of  $K^0_\Theta(H)$  are determined by homomorphisms $s:W(Y) \to
W(X)$ by the rule $\tilde s(\nu) = \nu s$ for every $\nu: W(X) \to
H$.  We have a contravariant functor
\[
\phi  : \Theta^0\to K^0_\Theta (H).
\]
\begin{prop} \mbox{[15]}
Functor $\phi$ determines duality  of categories if and only if
$Var (H) = \Theta$.
\end{prop}
\noindent {\bf Corollary} . If $Var(H_1) = Var (H_2) = \Theta$,
then the categories $K^0_\Theta(H_1)$ and $K^0_\Theta(H_2)$ are
isomorphic.

Proceed now to the category of algebraic sets $K_\Theta(H)$. Its
objects have the form $(X, A)$, where $A$ is an algebraic set in
the space $Hom(W(X), H)$.  The morphisms $[s]: (X, A) \to (Y, B)$
are defined by those $s: W(Y) \to W(X)$, for which $\tilde s(\nu)
\in B$ if $\nu \in A$.   Simultaneously, we have mappings
$[s]:A\to B$.

Let us define the category $C_\Theta(H)$.  Its objects are of the
form $W/T$, where $W\in 0b\Theta^0$ and $T$ is an $H$-closed
congruence in $W$.  Morphisms are homomorphisms of algebras.

It is proved that if $Var (H) = \Theta$ then the transitions $(X,
A) \to W(X)/A'$ and $W/T \to (X, T'_H)$ determine duality of the
categories $K_\Theta(H)$ and $C_\Theta(H)$.  In this case the
category $\Theta^0$ is a subcategory in $C_\Theta(H)$.

The category $K^0_\Theta (H) $ is always a subcategory in
$K_\Theta (H)$.

Regarding categories $K_\Theta$ and $C_\Theta$  see [21].
Correspondingly,  we have the categories $\tilde K_\Theta$ and
$\tilde C_\Theta$.

\subsection{Functor $Cl_H:\Theta^0 \to$ Set}

This functor corresponds to every algebra $H$ in $\Theta$.  By
definition, for every $W \in 0b\Theta^0$ the set $Cl_H(W)$ is the
set of all $H$-closed congruences $T$ in $W$.

Let now a morphism
\[
s:W(Y) \to W(X)
\]
be given in $\Theta^0$.  It corresponds a map
\[
Cl_H(s): Cl_H(W(X)) \to Cl_H(W(Y)),
\]
defined by the rule $Cl_H(s)(T) = s^{-1}T$.
 Here $T \in Cl_H(W(X)); s^{-1} T$ is a congruence in $W(Y)$,
defined by the rule $w(s^{-1}T) w'$ if and only if $w^sTw'^s,\; w,
w' \in W(Y)$.  The congruence $s^{-1} T$ is also $H$-closed.

This defines a contravariant   functor $Cl_H$, which plays an
important part in the sequel.

If $\Theta_1$ is a subvariety in $\Theta$, containing the algebra
$H$, then there is also $Cl_H: \Theta^0_1 \to Set$, see [20].

\section{Logically compact classes of algebras}

\subsection{Preliminary remarks}

We generalize the notion of $H$-closed congruence.  Let $\cx$ be
an arbitrary class of algebras in $\Theta$, $W = W(X) \in 0b\;
\Theta^0$.

\bigskip

\noindent{\bf Definition 1}.  The congruence $T$ in $W$ is called
$\cx$-closed, if $W/T \in SC(\cx)$.

It is clear that the intersection of $\cx$-closed congruences is
an $\mathcal \cx$-closed congruence as well, and, hence, for every
$T$ in $W$ one can consider $\cx$-closure, denoted by $T^\cx$. We
check directly that $T^\cx=\bigcap\limits_{H \in \cx} T^{''}_H$.

\bigskip

\noindent{\bf Proposition 3}.
{\it  The pair $(w_0, w'_0)$ belongs to $T^\cx$ if
and only if an infinitary quasiidentity
\begin{equation}\tag{$\ast$}
 \Big(
\bigwedge\limits_{(w, w') \in T} (w \equiv w')\Big) \Rightarrow
w_0 \equiv w'_0
\end{equation}
holds in the class $\cx$.}

\begin{proof}  Let the quasiidentity $(\ast)$ hold in $\cx$.
Algebra $W/T^\cx$ belongs to the class $LSC(\cx)$ by  definition.
By Theorem 1, the quasiidentity $(\ast)$ holds in this algebra.
Since the premise holds in it, so does the consequence. This means
that $(w_0, w'_0) \in T^\cx$.

Let now $(w_0, w'_0) \in T^\cx$. We need to verify that $(\ast)$
holds in $\cx$.  Take an arbitrary algebra $H\in \cx$.  We have
$SC(H)\subset SC(\cx)$.  Therefore, the algebra $W/T^{''}_H$
belongs to the class $SC(\cx)$. This gives the inclusion $T^\cx
\subset T^{''}_H$, and then $(w_0, w_0')\in T^{''}_H$.

According to  proposition 2.1 we may claim now that the
quasiidentity $(\ast)$ holds in the algebra $H$. Since $H$ is an
arbitrary algebra, the quasiidentity holds in $\cx$.
\end{proof}

\subsection{Logically compact classes}

We want to return to the problem: is
\[
\tilde q Var (\cx) = qVar (\cx)?
\]
\bigskip

\noindent {\bf Definition 2}. A class of algebras $\cx$ we call
logically compact $(q_w$-compact in [17] ), if  every infinitary
quasiidentity of this class is reduced in $\cx$ to an ordinary
(finitary) quasiidentity.
\medskip
This means that if $\ast$ holds in $\cx$ then there is a finite
subset $T_0$ in $T$ such that the same quasiidentity with the
equalities in $T_0$ holds in $\cx$.

\noindent {\bf Proposition 4}. {\it Class $\cx$ is logically
compact if and only if for every algebra $W=W(X) \in Ob \;
\Theta^0$ the union of a directed system of $\cx$-closed
congruences is also $\cx$-closed.}

\begin{proof}
System of congruences $T_\a, \a \in I$ is a directed system, if
for every $T_\a$ and $T_\beta$ there exists $T_\gamma$ with $T_\a,
T_\beta \subset T_\gamma $.  Union $T$ of the system of all $T_\a$
is a congruence.

Let now a class $\cx$ be logically compact and let all $T_\a$ be
$\cx$-closed. Let us show that $T$ is $\cx$-closed congruence.

Take a closure $T^\cx$ and let $(w_0, w'_0) \in T^\cx$.  Then we
have an infinitary quasiidentity $(\ast)$
\[
(\bigwedge_{(w, w') \in T}(w\equiv w')) \Rightarrow (w_0 \equiv
w'_0)\] in the class $\cx$.  Since $\cx$ is compact, then there
exists a finite subset $T_0$ in $T$, such that the quasiidentity
$(\ast)$ is equivalent in $\cx$  to the quasiidentity
\[
(\bigwedge_{(w, w') \in T_0} (w \equiv w')) \Rightarrow (w_0
\equiv w'_0).
\]
The set $T_0$ belongs to some $\cx$-closed $T_\a$.  Then $(w_0,
w'_0) \in T_\a \subset T$.  Therefore, $T^\cx = T$.

Let now the condition about directed systems of congruences hold.
Let us prove that $\cx$ is a logically compact class.

Take an infinitary quasiidentity
\[\tag{$\ast$}
(\bigwedge_{(w, w')\in T} (w\equiv w')) \Rightarrow (w_0 \equiv
w'_0)
\] and let it hold in $\cx$.  Let $(\ast)$  be defined over
$W=W(X)$.  Consider various finite subsets $T_\a$ of the set $T$,
and for every $T_\a$ pass to $T^\cx_\a$.  All $T^\cx_\a$
constitute a directed system of $\cx$-closed congruences in the
algebra $W$. Let $T_1$ be the union of all $T^\cx_\a$.  We have:
$T^\cx_1 = T_1, \; T^\cx \subset T_1$. Since quasiidentity
$(\ast)$ holds in the class $\cx$, then $(w_0, w'_0) \in
T^\cx\subset T_1$. Hence, $(w_0, w'_0)$ is contained in some
$T^\cx_\a$.  This means, that the finitary quasiidentity
\[\tag{$\ast\ast$} (\bigwedge_{(w, w') \in T_\a} (w\equiv w')) \Rightarrow (w_0
\equiv w'_0)
\]
holds in the class $\cx$. The initial infinitary quasiidentity
$(\ast)$ is reduced to $(\ast\ast)$.  Class $\cx$ is logically
compact.
\end{proof}

\begin{thm} \mbox{[17]}
The equality $\tilde q Var (\cx) = q Var (\cx)$ holds if and only
if the class $\cx$ is logically compact.
\end{thm}
\begin{proof}
This theorem has been proved for groups in [17] . It was noted
there that the proof is valid for any $\Theta$. We present the
proof for an arbitrary variety of algebras $\Theta$, which is
slightly different from the proof from [17].  Note first of all
that it follows from the definitions that if $\cx$ is a logically
compact class, then
\[
\tilde q Var (\cx) = q Var (\cx).
\]
Now let this equality hold true.  Let us prove that the class
$\cx$ is logically compact.

Take an arbitrary algebra $W=W(X)$ and prove that if $T_\a, \; \a
\in I$ is a directed system of $\cx$-closed congruences in $W$,
and $T$ is the union of this system, then the congruence $T$ is
also $\cx$-closed.  This implies that the  class $\cx$ is
logically compact.

Every algebra $W/T_\a$ belongs to the class $SC(\cx)$.  We need to
check that $W/T$ belongs to this class as well.  Since $LSC(\cx) =
q Var (\cx)$ and algebra $W$ is finitely generated, it is enough
to check that all the quasiidentities of the class $\cx$ hold in
the algebra $W/T$.  Let
\begin{equation}\tag{$\ast\ast\ast$}
w_1\equiv w'_1\wedge\dots \wedge w_n\equiv w'_n \Rightarrow w_0
\equiv w'_0
\end{equation} be one of such quasiidentities, written in the
algebra $W(Y)$.

Consider an arbitrary homomorphism
\[
\mu: W(Y) \to W(X)/T
\]
and associate to it a commutative diagram with $\mu$:

$$
\CD
W(Y) @>\mu_0>> W(X)\\
@[2]/SE/\mu // @. @VV\nu V\\
@. W(X)/T\\
\endCD
$$
Here $\nu $ is the natural homomorphism.  Besides, for every $\a
\in I$ we have natural homomorphisms $\nu_\a: W(X) \to W(X)/T_\a$.
Assume that $w_i^\mu = w^{'\mu}_i, w^{\mu_0\nu}_i= w_i^{'\mu_0
\nu}$ holds for every $i_1, \dots, i_n$.  We can choose $\a \in I$
such that $w_i^{\mu_0\nu_\a}=w_i^{'\mu_0 \nu_\a}$ holds also for
every $i = 1, \dots, n$.  We proceed from the homomorphism $\nu_\a
\mu_0: W(Y) \to W(X)/T_\a$.  Since the algebra $W(X)/T_\a$ belongs
to the class $LSC(\cx)$, the quasi-identity $(\ast\ast\ast)$ holds
in it.  Hence, $w_0^{\mu_0\nu_\a} = w_0^{'\mu_0\nu_\a}$.  Then
$w_0^{\mu_0 \nu} = w_0^{'\mu_0\nu}, w^\mu_0 = w_0^{'\mu}$.  This
means that the quasiidentity $(\ast\ast\ast)$ holds in $W(X)/T$
and the congruence $T$ is $\cx$-closed.  Hence, the class $\cx$ is
logically compact.

We now return to geometric notions.
\end{proof}
\section{Geometrically equivalent algebras}

\subsection{Definition}
Algebras $H_1$ and $H_2$ from
$\Theta$ are called geometrically equivalent if for every $W=W(X)
\in Ob\; \Theta^0$ and every $T$ in $W$, we have
\[
T^{''}_{H_1} = T^{''}_{H_2}.
\]
This means also that $Cl_{H_1} = Cl_{H_2}$.   It is clear that if
the algebras $H_1$ and $H_2$ are geometrically equivalent, then
the categories $C_\Theta(H_1)$ and $C_\Theta(H_2)$ coincide.
Correspondingly,  the categories $K_\Theta(H_1)$ and
$K_\Theta(H_2)$ are isomorphic.
\begin{thm} \mbox{[24]}
Algebras $H_1$ and $H_2$ are geometrically equivalent if and only
if
\[ LSC(H_1) = LSC(H_2).
\]
Hence, geometrical equivalence of algebras means also that
\[
\tilde q Var (H_1) = \tilde q Var (H_2),
\]
i.e., $H_1$ and $H_2$ have the same infinitary quasiidentities.
\end{thm}
\medskip

\noindent {\bf Corollary}. {\it If $H_1$ and $H_2$ are
geometrically equivalent, then $q Var (H_1) = q Var(H_2)$ and $Var
(H_1) = Var (H_2)$.  The problem whether $q Var (H_1) = q Var
(H_2)$ implies geometrical equivalence of $H_1$ and $H_2$ has the
negative answer  (Theorem 7 in the sequel).}

\subsection{Twisted and almost geometrically equivalent algebras}

Let algebra $H$ belong to $Ass-P$ or $Lie-P$ and $\sigma \in Aut
(P)$.  Define a new algebra $H^\sigma$.   In $H^\sigma$ the
multiplication on a scalar $\circ$ is defined through the
multiplication in $H$ by the rule:
\[
\lambda\circ a=\lambda^{\sigma^{-1}}\cdot a, \;\; \;   \lambda \in
P, \;\; \;  a \in H.
\]
I.e., $ \lambda a = \lambda^\sigma \circ a.$ We say that the
algebra $H^\sigma$ is $\sigma$-twisted in respect to $H$. This is
also an associative or Lie algebra. Besides, note that the
identical map $H\to H^\sigma$ is a semiisomorphism of algebras.
\bigskip

\noindent{\bf Definition 3}.
Algebras
 $H_1$ and $H_2$ are called twisted geometrically equivalent if
 $H^\sigma_1$ and $H_2$ are geometrically  equivalent for some
 $\sigma$.

 Let, further, $H$ be an associative algebra.  Denote by $H^*$ an
 opposite algebra
 \[
 a\circ b = b a.
 \]
 An identical transition $H\to H^*$ here is an antiisomorphism of
 algebras.
\bigskip

\noindent{\bf Definition 4}.
 Associative  algebras $H_1$ and $H_2$ are called almost
 geometrically equivalent if they are twisted geometrically
 equivalent or  $ (H_1^\sigma)^*$ and $H_2$ are geometrically
 equivalent for some $\sigma \in Aut (P)$.

 \section{Main results}
 \subsection{$\Theta=Com-P$}

 Let the field  $P$ be infinite.  In this case we have $Var(H) = \Theta$
 for every $H\in\Theta$.

 Besides, we will see that algebras $H_1$ and $H_2$ in $\Theta$
 are geometrically equivalent if and only if they have the same
 quasiidentities.  Hence, geometrical equivalence in $Com-P$ is, in
 fact, logical equivalence in the logic of quasiidentities.
\begin{thm}
Let $H_1$ and $H_2$ be algebras in $\Theta=Com-P$.  Then the
following three conditions are equivalent:
\begin{enumerate}
\item[1.] Categories $K_\Theta(H_1)$ and $K_\Theta(H_2)$ are
isomorphic.
\item[2.] Categories $K_\Theta (H_1)$ and $K_\Theta(H_2)$ are
equivalent.
\item[3.] $H_1$ and $H_2$ are twisted geometrically equivalent.
\end{enumerate}
\subsection{$\Theta=Lie-P$}
\end{thm}
\begin{thm}
Let $Var(H_1) = Var(H_2) = \Theta$. Then the following conditions
are equivalent:
\begin{enumerate}
\item[1.] $K_\Theta(H_1)$ and $K_\Theta(H_2)$ are isomorphic
\item[2.] $K_\Theta(H_1)$ and $K_\Theta(H_2)$ are equivalent.
\item[3.] $H_1$ and $H_2$ are twisted geometrically equivalent.
\end{enumerate}
\end{thm}
\subsection{$\Theta=Ass-P$}

\noindent{\bf Conjecture 1:}

 Let $Var(H_1)= Var (H_2) = \Theta$.  Then the following
conditions are equivalent:
\begin{enumerate}
\item[1.] $K_\Theta(H_1)$ and $K_\Theta(H_2)$ are isomorphic.
\item[2.] $K_\Theta(H_1)$ and $K_\Theta(H_2)$ are equivalent
\item[3.] Algebras $H_1$ and $H_2$ are almost geometrically
equivalent.
\end{enumerate}

There will be some observations in the sequel in favor of this
conjecture.

\subsection{Correctness}

In fact, we use here  special correct isomorphisms and correct
equivalences.  As we will see, these notions are natural and
reflect well the idea of coincidence of geometries.

Roughly speaking, correctness means correlation with the category
of affine spaces and with inclusions of algebraic sets.

More precisely, let an isomorphism $F: K_\Theta(H_1) \to K_\Theta
(H_2)$ be given.  Then an isomorphism $\Phi: C_\Theta (H_1) \to
C_\Theta (H_2)$ corresponds to it.  Correctness of $F$ assumes,
that
\begin{enumerate}
\item[1)] $\Phi$ induces an automorphism $\Phi_0 = \varphi:
\Theta^0\to \Theta^0$ and $F$ induces an isomorphism $F_0:
K^0_\Theta (H_1) \to K^0_\Theta (H_2)$.
\item[2)] Let $(X, A_1)$ and $(X, A_2)$ be  two objects of
$K_\Theta(H_1), \; A_1 \subset A_2$, and let $F(X, A_1)= (Y_1,
B_1), \; F(X, A_2) =(Y_2, B_2)$.  Then $Y_1=Y_2=Y, \;$ and $
B_1\subset B_2$.
\end{enumerate}

It follows from this definition, that a correct isomorphism $F$ is
well coordinated with the lattices of algebraic sets.

The correctness of an equivalence of categories is defined in the
same spirit.

\subsection{Program of further considerations and proofs}
Our plan is as follows:

\begin{enumerate}
\item[1.] Investigate the notion of geometrical equivalence
in more details.
\item[2.] Generalize this notion and consider the notions of
geometrically similar algebras and geometrically compatible
algebras.
\item[3.] Prove the universal theorems (for arbitrary variety $\Theta$)
 about isomorphism and equivalence of categories of algebraic
sets.  We use here the  notions of geometric similarity and
geometric compatibility.
\item[4.] In order to apply these universal theorems to specific
$\Theta$ we need information about automorphisms of the category
$\Theta^0$ of free in $\Theta$ algebras.
\item[5.] Apply these
four steps to the cases $Com-P$, $Ass-P$, $Lie-P$.
\end{enumerate}

\section{Geometrically equivalent algebras (continuation)}
\subsection{Geometrically noetherian algebras}

\noindent {\bf Definition 5}. An  algebra $H\in\Theta$ is
geometrically noetherian if for an arbitrary $W$ and $T$ in $W$
there exists a finite $T_0\subset T$ such that
 $$T^{''}_H = (T_0)^{''}_H.$$

\medskip

An algebra $H$ is geometrically noetherian if and only if in every
$W=W(X)$ the ascending chain condition for $H$-closed congruences
holds. The same holds for descending chain condition for algebraic
sets in $Hom(W(X), H)$.
\bigskip

\noindent{\bf Definition 6}.  The variety $\Theta$ is noetherian
if every $W=W(X) \in Ob\Theta^0$ is noetherian (in respect to
congruences).

\medskip

If $\Theta$ is noetherian then  an arbitrary algebra $H\in\Theta$
is geometrically noetherian.

Examples:
\begin{enumerate}
\item[1.] The variety $Com-P$ is noetherian
\item[2.] Free group of finite rank is geometrically noetherian
(Guba [8])
\item[3.] Associative and Lie algebras of finite dimension are
geometrically noetherian
\item[4.] The variety of nilpotent groups $N_c$ is noetherian
\item[5.] All noetherian subvarieties in $Ass-P$ are described in [1]
\end{enumerate}

\noindent {\bf Problem 4}

What are all noetherian subvarieties in the variety of all groups?

\subsection{Logically noetherian algebras}

\bigskip

\noindent{\bf Definition 7}. An algebra $H$ is called logically
noetherian if the class $\cx = \{ H\}$ is logically compact.
\medskip

Every geometrically noetherian algebra is logically noetherian.

\begin{thm} Let $H_1$ and $H_2$ be logically noetherian algebras.
 They are
geometrically equivalent if and only if $q Var (H_1) = q Var
(H_2)$.
\end{thm}
It follows from  Theorem 2 that the equality $LSC(H) = qVar (H)$
holds if and only if $H$ is logically noetherian.  This together
with the presentation for $qVar (H)$ imply
\begin{thm} {\rm{\bf ([17] )}}
If $H\in\Theta$ is not logically noetherian, then there exists an
ultrapower $H'$ of $H$ such that the algebras $H$ and $H'$ are not
geometrically equivalent.

 However, these algebras have the same elementary theories and, in
particular, the same quasiidentities.
\end{thm}

In [23] one can find examples of not logically noetherian groups
and associative algebras.  The results from Gobel-Shelah [5] and
Lichtman-Passman [11] are used in the proofs.

\noindent {\bf Problem 5}

Build examples of not logically noetherian Lie algebras.

\noindent {\bf Problem 6}

Let $W =W(X)$ be a free associative or free Lie algebra, $X$ is
finite.  Whether $W$ is not geometrically noetherian, but
logically noetherian.

\section{Geometrically similar algebras}

\subsection{Some information from category theory}

Recall first, that $s:\vp_1\to \vp_2$  is an isomorphism of two
covariant functors $\vp_1, \vp_2: C_1\to C_2$
 if to every
$A\in Ob \; C_1$ it corresponds the isomorphism
\[
s_A: \vp_1(A) \to \vp_2(A)
\]
in $C_2$ and for $\nu: A\to B$ in $C_1$  there is a commutative
diagram
$$
\CD
\vp_1(A) @>S_A>> \vp_2(A) \\
@V\varphi_1(\nu)VV @VV\vp_2(\nu)V\\
\varphi_1(B)@>S_B>> \varphi_2(B)\\
\endCD
$$

For contravariant functors $\vp_1$ and $\vp_2$ the corresponding
diagram looks as follows
$$
\CD
\vp_1(B) @>S_B>> \vp_2(B) \\
@V\vp_1(\nu) VV @VV \vp_2(\nu)V\\
\vp_1(A) @>S_A>> \vp_2(A) \\
\endCD
$$

Denote the relation of isomorphism  by $\approx$.

An endomorphism $\vp$ of the given category $C$ (covariant
endofunctor) we call an inner endomorphism, if there exists an
isomorphism $s:1_C\to \vp$.  For every $A\in Ob\;  C$ we have an
isomorphism $s_A: A  \to \vp(A)$ and for $\nu: A\to B$ a
commutative diagram
$$
\CD
A @>S_A>> \vp(A) \\
@V\nu VV @VV \vp(\nu)V\\
B @>S_B>> \vp(B)\\
\endCD
$$
holds.  Now,
\[
\vp(\nu)=s_B\; \nu s_A^{-1}: \vp (A) \to \vp(B).
\]
This motivates the name ``inner".  In particular, one can speak of
inner automorphisms of the given category.

It is easy to show that if $C$ is a monoid, considered as a
category, then inner automorphisms of this category are exactly
inner automorphisms of the monoid.

Recall now the definition of equivalence of two categories.

Consider a pair of functors:
\[
\vp: C_1\to C_2,\quad \psi:C_2\to C_1.
\]
The pair $(\vp, \psi)$ determines equivalence of categories if
$\psi\vp\approx 1_{C_1}, \; \vp \psi\approx 1_{C_2}$. Here
$\psi\vp$ and $\vp\psi$  are inner endomorphisms of the
corresponding categories $C_1$ and $C_2$.

If $C_1 = C_2 = C$, then $(\vp, \psi)$ is called autoequivalence
of the category $C$.

An autoequivalence $(\vp, \psi)$ we call an inner autoequivalence,
if $\vp$ and $\psi$ are inner.  In particular, if $\vp$  is inner,
then the pair $(\vp, \psi)$ is an inner autoequivalence.

\subsection{Definition of similarity}

We assume that the condition $Var(H_1) = Var(H_2) = \Theta$ holds
true.  This condition always holds in the classical situation.
Recall that the correctness of the isomorphism $F:K_\Theta(H_1)
\to K_\Theta(H_2)$ assumes that an automorphism $\vp:
\Theta^0\to\Theta^0$ corresponds to $F$.  Now we proceed from such
an automorphism and consider a diagram of functors:
$$
\CD
\Theta^0 @[2]>\vp>> \Theta^0 \\
 @[2]/SE/Cl_{H_1} //@.@.\;    @/SW//Cl_{H_2} / \\
 @. Set \\
\endCD
$$

Commutativity of the diagram means that to the automorphism $\vp$
it corresponds the transition
\[
\a (\vp): Cl_{H_1} \to Cl_{H_2} \cdot \vp,
\]
with the following properties:
\begin{enumerate}
\item[1.] To every $W=W(X) \in Ob\; \Theta^0$ it corresponds  the
bijection
\[
\a(\vp)_W: Cl_{H_1}(W) \to Cl_{H_2} (\vp(W))
\]
\item[2.]  The function $\a(\vp)$ should be compatible
with the automorphism $\vp$.
\end{enumerate}
Let us explain the condition of $\vp$ and $\a(\vp)$ compatibility.
Let $W_1$ and $W_2$ be objects of $\Theta^0$,
\[ T\in C1_{H_1} (W_2),\; \; \,  T^*=\a (\vp)_{W_2} (T) \in C1_{H_2}
\vp(W)\] and let $\mu_T: W_2 \to W_2/T$ and $\mu_{T^*}: \vp (W_2)
\to \vp(W_2) /T^*$ be natural homomorphisms.  Then for any $s_1,
s_2: W_1\to W_2$ it should hold: the equality $\mu_T s_1=
\mu_Ts_2$ fulfills if and only if $\mu_{T^*} \vp (s_1) = \mu_{T^*}
\vp (s_2))$.
\bigskip

\noindent{\bf Definition 8}. Algebras $H_1$ and $H_2$ are
geometrically similar if for some $\vp: \Theta^0 \to \Theta^0$ the
above conditions  hold.

We say that the automorphism $\vp$ determines similarity of
algebras.  Properties of this $\vp$ determine properties of
similarity.  In some cases similarity is reduced to geometrical
equivalence, or to twisted equivalence, or to almost equivalence.

For the identical $\vp$ we have geometrical equivalence.  Here
$\a(\vp)$ determines the equality $Cl_{H_1} = Cl_{H_2}$.

\subsection{Corollary from the definition}

Note first of all the following theorem:
\begin{thm}\mbox{[23]}
The transition $\a (\vp): Cl_{H_1} \to Cl_{H_2} \cdot \vp$ is an
isomorphism of functors.
\end{thm}
 Let us study the structure of
this isomorphism.  Let $W = W(X) \in Ob\; \Theta^0$ and $T$ be a
congruence in $W$.  The relation $\rho= \rho(T) = \rho_W(T)$ in
the semigroup $End\,  W$ is determined  by the rule
\[
\nu\rho\nu'\Leftrightarrow \mu_T\nu = \mu_T\nu', \; \, \; \nu, \nu' \in
End \, W.\]

Let, further, $\rho $ be an arbitrary binary relation in $End \,
W$. Define the relation $\tau = \tau(\rho)= \tau_W(\rho)$ in $W$
by the rule:
\[
w_1\tau w_2\Leftrightarrow \exists w, \;  \nu, \nu' \ | \
 w_1 = w^\nu,\; \; \,  w_2 = w^{\nu'}, \; \nu\rho\nu'.
 \]
 If $T$ is a congruence, then $\tau_W(\rho_W(T))=T$.  It follows
 from the definitions, that
 \[
 \a(\vp)_W(T) = \tau_{\vp (W)} (\vp(\rho_WT)).\]
Here $\vp(\rho)$ is a relation in $End \vp (W)$ defined by the
rule: if $\mu, \mu' \in End \vp (W)$, then $\mu\vp( \rho) \mu'$ if
and only if there exist $\nu$ and $\nu' \in End \, W$ with $\vp
(\nu) = \mu, \vp (\nu') = \mu'$ and $\nu\rho\nu'$.  This gives
rise to the proof that the transition
\[
\a (\vp)_W: Cl_{H_1} (W) \to Cl_{H_2} (\vp(W))\] is an isomorphism
of lattices of algebraic sets in $Hom(W, H_1)$ and $Hom(\vp(W),
H_2)$.

\subsection{The main theorem}
\begin{thm} \mbox{[23]}  Categories $K_\Theta (H_1)$ and $K_\Theta(H_2)$ are
correctly isomorphic if and only if the algebras  $H_1$ and $H_2$
are geometrically similar.
\end{thm}
This theorem, as well as the similar theorem on correct
equivalence of categories, is used in special cases $Com-P, Ass-P$
and $Lie-P$.
\section{Geometrical compatibility of algebras}
\subsection{Definition}
As earlier, we consider the diagrams of functors

$$
\CD
\Theta^0 @[2]>\vp>> \Theta^0 \\
 @[2]/SE/Cl_{H_1} //@.@.\;    @/SW//Cl_{H_2} / \\
 @. Set \\
\endCD
$$

$$
\CD
\Theta^0 @[2]<\psi<< \Theta^0 \\
 @[2]/SE/Cl_{H_1} //@.@.\;    @/SW//Cl_{H_2} / \\
 @. Set \\
\endCD
$$



Here the pair $(\vp,\psi):\Theta^0\to \Theta^0$ determines
autoequivalence  of the category $\Theta^0$. Suppose that the
transitions
\begin{align*}
\a(\vp): \; &Cl_{H_1} \to Cl_{H_2}\;  \vp,\\
\a(\psi):\; &Cl_{H_2} \to Cl_{H_1} \; \psi
\end{align*}
are given.  Then for every $W\in Ob\; \Theta^0$ we have the
mappings
\begin{align*}
\a(\vp)_W: \; &Cl_{H_1}(W) \to Cl_{H_2} (\vp(W)),\\
\a(\psi)_W:\; &Cl_{H_2}(W)  \to Cl_{H_1} (\psi(W)).
\end{align*}
We assume also, that these mappings are compatible with the
initial autoequivalence $(\vp, \psi)$ like it was in the
definition of geometrical similarity.
\bigskip

\noindent{\bf Definition 9}. Algebras $H_1$ and $H_2$ are
geometrically compatible by the autoequivalence $(\vp,\psi)$ if
there are $\a(\vp)$ and $\a(\psi)$ for $(\vp, \psi)$, which
satisfy the compatibility conditions.

\subsection{Corollaries from the definition}

First of all note that the transitions  $\a(\vp)_W$ can be
presented in the form
\begin{align*}
\a(\vp)_W(T) =  &\tau_{\vp(W)}\;  \vp(\rho_W(T))\\
\a(\psi)_W(T)= &\tau_{\psi(W)} \; \psi(\rho_W(T))
\end{align*}
.

We deduce from the definitions the following
\begin{thm}
The transitions
\begin{align*}
\a(\vp): \; &Cl_{H_1} \to Cl_{H_2}\;  \vp,\\
\a(\psi):\; &Cl_{H_2} \to Cl_{H_1} \; \psi
\end{align*}
are natural transformations (morphisms) of functors.
\end{thm}
\begin{proof}
It is enough to study the transition $\a(\vp)_W: Cl_{H_1} \to
Cl_{H_2} \vp$.  Proceed from the morphism $s:W_1\to W_2$ and
consider the diagram
$$
\CD
Cl_{H_1}(W_2) @>\a(\vp)_{W_2}>> Cl_{H_2}(\vp(W_2)) \\
@V Cl_{H_1} (s)VV @VV Cl_{H_2} \vp(s) V\\
Cl_{H_1}(W_1) @>\a(\vp)_{W_1}>>
Cl_{H_2} (\vp(W_1))\\
\endCD
$$
Check the commutativity of this diagram:
\[Cl_{H_2} (\vp(s)) \a
(\vp)_{W_2} = \a(\vp)_{W_1} Cl_{H_1} (s)\] Apply both parts to
$T\in Cl_{H_1} (W_2)$.  We have $Cl_{H_1}(s) (T) = s^{-1} T$.
Denote $\a(\vp)_{W_2}(T) = T^* \in Cl_{H_2} (\vp(W_2))$. Then
$Cl_{H_2} \vp(s)(T^*)= \vp(s)^{-1} (T^*)$. So we need to verify
that
\begin{align*}
&\vp(s)^{-1} T^* =\a (\vp)_{W_1} (s^{-1}T),\\
&\vp(s)^{-1} \a(\vp)_{W_2} (T) = \a (\vp)_{W_1} (s^{-1} T).
\end{align*}
Both parts belong to $Cl_{H_2} (\vp(W_1))$.  Take elements $w_1,
w_2$ in $\vp(W_1)$.  Given $w_1(\a(\vp)_{W_1} (s^{-1} T)) w_2$, we
have
\[
\a (\vp)_{W_1} (s^{-1} T) = \tau_{\vp (W_1)} \vp (\rho_{W_1}
(s^{-1} T)).\] It follows from the definition of the function
$\tau$, that there exist $w_0 \in \vp(W_1)$, $ \mu_1, \mu_2 \in
End (\vp(W_1))$ such that $w_1= w_0^{\mu_1}, w_2= w_0^{\mu_2},
\mu_1\vp(\rho_{W_1} (s^{-1} T)) \mu_2$.  We use the univalence
property of the functor $\vp$ [16].  According to this property,
there exist unique $\nu_1, \nu_2 \in End (W_1)$ with $\vp(\nu_1) =
\mu_1, \vp(\nu_2) = \mu_2$.   Now the condition $\vp(\nu_1) \vp
(\rho_{W_1} (s^{-1} T)) \vp (\nu_2)$ means that $\nu_1\rho_{W_1}
(s^{-1} T) \nu_2$ holds.  For every $w \in W_1$ we have $w^{\nu_1}
(s^{-1} T) w^{\nu_2}$.  Then $w^{\nu_1 s} Tw^{\nu_2 s}$ which is
equivalent to $\mu_T(s\nu_1) = \mu_T(s\nu_2)$ for a natural
homomorphism $\mu_T:W_1\to W_2/T$.

Applying the condition of compatibility of the function $\a(\vp)$
and the functor $\vp$, we get
\[
\mu_{T^*} \vp(s\nu_1) = \mu_{T^*} \vp(s\nu_2).\] The latter means
that for every $w\in \vp(W_1)$ there hold
\begin{align*}
&w^{\vp(\nu_1 s) } T^* w^{\vp(\nu_2s)},\\
&w^{\vp(\nu_1)\vp(s)}T^*w^{\vp(\nu_2)\vp(s)},\\
&w^{\vp(\nu_1)} \vp(s)^{-1} T^*w^{\vp(\nu_2)}
\end{align*}
Apply this to the initial $w_0, \mu_1, \mu_2$:
\begin{align*}
&w_0^{\mu_1} \vp(s)^{-1} T^* w_0^{\mu_2},\\
&w_1(\vp(s)^{-1} \a(\vp)_{W_2} (T)) w_2.
\end{align*}
We have checked
\[
\a(\vp)_{W_1} (s^{-1} T) \subset \vp(s)^{-1} \a (\vp)_{W_2} (T).
\]
Check now the opposite inclusion.  Let $w_1(\vp(s)^{-1} \a
(\vp)_{W_2} (T)) w_2$.  Take $T_0 = \vp (s)^{-1} \a (\vp)_{W_2}
(T) \in Cl_{H_2} \vp (W_1)$.  We have $T_0 = \tau_{\vp(W_1)}
\rho_{\vp (W_1)}(T_0)$.  Here $w_1 T_0 w_2$ means that in
$\vp(W_1)$ there exists $w_0$ and in $End \vp(W_1)$ there exist
$\mu_1, \mu_2$ such that
\[
w_0^{\mu_1} = w_1,\ w_0^{\mu_2}=w_2, \ \mu_1 \rho_{\vp(W_1)}(
T_0)\mu_2.\] For every  $w\in \vp (W_1)$ we have $w^{\mu_1} T_0
w^{\mu_2}$; $w^{\mu_1\vp(s)} \a(\vp)_{W_2} (T) w^{\mu_2\vp(s)}$.
Taking into account univalencity of the function $\vp$, we get
$\nu_1, \nu_2\in End\,  W_1 $ with $\vp (\nu_1) = \mu_1,
\vp(\nu_2) = \mu_2$.  This gives $w^{\vp(\nu_1 s)} (\a(\vp)_{W_2}
(T)) w^{\vp(\nu_1s)}$.  Taking $T^* = \a (\vp)_{W_2}(T)$, we come
to
\[
\mu_{T^*} (\vp(\nu_1 s)) = \mu_{T^*} (\vp(\nu_2s)).
\]
Compatibility condition for $\vp $ and $\a (\vp)$ implies  $\mu_T
\nu_1s = \mu_T \nu_2 s$.  For every $w \in \vp (W_1)$ it holds
\begin{align*}
w^{\nu_1 s} &Tw^{\nu_2s};\ \  w^{\nu_1} (s^{-1} T) w^{\nu_2};\\
&\nu_1\rho_{W_1} (s^{-1} T) \nu_2.
\end{align*}
Take now $T_1 = s^{-1} T$ and apply the condition of compatibility
of $\vp$ and $\a (\vp)$ in the case $W_2 = W_1$:
\[
\mu_{T_1} \nu_1 = \mu_{T_1} \nu_2\Rightarrow \mu_{T^*_1}
\vp(\nu_1)= \mu_{T^*} \vp(\nu_2).
\]
The condition $w^{\vp(\nu_1)} \a (\vp)_{W_1} (s^{-1} T)
w^{\vp(\nu_1)}$ holds for every $w \in \vp(W_1)$.  Applying this
to the initial $w=w_0, \mu_1$ and $\mu_2$, we get
\begin{align*}
&w^{\mu_1}_0 \a (\vp)_{W_1} (s^{-1} T) w_0^{\mu_2};\\
&w_1(\a(\vp)_{W_1} (s^{-1}T)) w_2.
\end{align*}
We checked the opposite inclusion.  The theorem is proved.
\end{proof}
\subsection{The main theorem}
\begin{thm}\mbox{[23]}
Categories $K_\Theta(H_1) $ and $K_\Theta(H_2)$ are correctly
equivalent if and only if algebras $H_1$ and $H_2$ are
geometrically compatible.
\end{thm}
Keeping in mind further applications of  theorems 9 and 11, let us
pass to automorphisms and autoequivalences of categories.
\section{Automorphisms and autoequivalences of categories of free
algebras of varieties. Applications}

\subsection{Semigroups $ End C$ and $End^0 C$}

Let $C$ be an arbitrary small category.  Denote by $End(C)$ a
semigroup of all endomorphisms (covariant endofunctors) of this
category.  It can be verified that the relation of isomorphism of
functors $\approx$ is a congruence of the semigroup $End(C)$.
Denote $End^0 (C) = End (C)/\approx$.  We have a natural
homomorphism $\delta: End(C) \to End^0(C)$.  The group $Aut(C)$ is
the group of invertible elements in $End(C)$, and $Aut^0(C)$ is
the group of invertible elements   in $End^0(C)$.  We have a
homomorphism $\delta: Aut(C) \to Aut^0(C)$.  The kernel
$Ker(\delta) = Inn (C)$ consists of inner automorphisms.
\subsection{Categories of $\Theta^0$ type}

\begin{thm}\mbox{[26]}
If the pair $(\vp, \psi)$ is an autoequivalence of the category
$\Theta^0$, then $\vp=\vp_0 \zeta, \; \psi=\zeta^{-1} \psi_0$,
where $\zeta$ is an automorphism of $\Theta^0$ and $(\vp_0,
\psi_0)$ is an inner autoequivalence of the category $\Theta^0$.
From this follows that the homomorphism $\delta: Aut (\Theta^0)\to
Aut^0(\Theta^0)$ is a surjection.
\end{thm}
\subsection{Special categories $\Theta$}

Let at the beginning $\Theta=Ass-P$ or $\Theta=Lie-P$. In these
cases we can consider semiinner automorphisms and
autoequivalences.

An automorphism $\vp:\Theta^0 \to \Theta^0$ we call semiinner if
it is semiisomorphic to an identity functor.  This means that
there is a semiisomorphism $(\sigma, s):1_{\Theta^0} \to \vp$,
where $\sigma \in P$ and for every $W\in Ob\;(\Theta^0)$ there is
$\sigma$-semiisomorphism $(\sigma, s_W): W \to \vp(W)$.  Besides
that, $\vp(\nu)= s_{W_2} \nu s^{-1}_{W_1}: \vp (W_1) \to \vp
(W_2)$ for $\nu: W_1 \to W_2$.

Let us define a mirror automorphism of the category $\Theta^0$ for
$\Theta=Ass-P$.

Let $W=W(X)$ be a free associative algebra over a field $P$, i.e.,
the algebra of noncommutative polynomials and $X$ a finite set.
Let $S(X)$ be a free semigroup  over $X$, and $W(X) = PS(X)$ a
semigroup algebra.  Every element $w \in W(X)$ has the form
\[
w=\lambda_0+\lambda_1 u_1 + \ldots + \lambda_k u_k, \; \; \;
\lambda_i\in P,
\]
all $u$ lie in $S(X)$.  Let now $u = x_{i_i}\ldots x_{i_n}$.  Take
$\bar u = x_{i_n} \ldots x_{i_1}$.  For $w$ take
\[
\bar w=\lambda_0+\lambda_1\bar u_1 + \ldots + \lambda_k\bar u_k.
\]
Define now a mirror automorphism $\eta$ of the category
$\Theta^0$, $\Theta = Ass-P$.  For  every $W\in Ob\;\Theta^0$ we
have $\eta(W) = W$.  Objects are not changed.

Let $\nu:W(X) \to W(Y)$ be given.  Define $\eta(\nu): W(X) \to
W(Y)$, setting $\eta(\nu) (x)=\overline{\nu(x)}$ for every $x\in
X$.

The following theorem takes place:
\begin{thm} \mbox{[14]}
For particular varieties $\Theta$ we have
\begin{enumerate}
\item[1.] $\Theta= Grp$, all automorphisms of the category
$\Theta^0$ are inner
\item[2.] $\Theta = $ semigroups, $Inn(\Theta^0)$ has index 2
in $Aut (\Theta^0)$.
\item[3.] $\Theta=Grp-F$, all automorphisms are semiinner.
\item[4.] $\Theta = Com-P$, all automorphisms are semiinner.
\item[5.] $\Theta=Lie-P$, all automorphisms are semiinner.
\item[6.] $\Theta=Mod-K, \; \; K$ is left noetherian, all
automorphisms of the category $\Theta^0$ are semiinner.
\end{enumerate}
\end{thm}

\noindent{\bf Conjecture 2}. (Special case $\Theta=Ass-P$.)

{\it{ All automorphisms of the category $\Theta^0$ are either
semiinner, or of the type $\vp_0 \eta$, where $\vp_0$ is semiinner
and $\eta $ is mirror.}}

The corresponding reduction theorem [15] allows to reduce this
case to the study of the group $Aut (End\, W(x, y))$.  Here $W(x,
y)$ is the free associative algebra with two variables. Positive
answer on this conjecture allows to answer positively on the main
conjecture when the geometries for algebras in $Ass-P$ are the
same.

The proof of the principal theorems for $Com-P, \; Ass-P$ and
$Lie-P$ is based on the following theorem:
\begin{thm}
Let similarity algebras $H_1$ and $H_2$ be determined by  an
automorphism $\vp$.  Then
\begin{enumerate}
\item[1.] If $\vp$ is inner, then $H_1$ and $H_2$ are
geometrically equivalent.
\item[2.] If $\vp$ is semiinner, then $H_1$ and $H_2$ are twisted
equivalent.
\item[3.] If $\vp = \vp_0\eta$, \ $\vp_0$ is semiinner, then $H_1$
and $H_2$ are almost geometrically equivalent.
\end{enumerate}
\end{thm}

Analogous theorem takes place for the relation of compatibility of
$H_1$ and $H_2$, determined by an autoequivalence $(\vp, \psi)$ of
the category $\Theta^0$.

\end{document}